\documentclass{amsart}

\usepackage{amssymb,latexsym}

\newtheorem{theorem}{Theorem}
\newcommand{\bt}{\begin{theorem}}
\newcommand{\et}{\end{theorem}}

\newtheorem{lemma}{Lemma}
\newcommand{\bl}{\begin{lemma}}
\newcommand{\el}{\end{lemma}}
\newtheorem{corollary}{Corollary}
\newcommand{\bc}{\begin{corollary}}
\newcommand{\ec}{\end{corollary}}

\newcommand{\beq}{\begin{equation}}
\newcommand{\eeq}{\end{equation}}
\newcommand{\benum}{\begin{enumerate}}
\newcommand{\eenum}{\end{enumerate}}

\newcommand{\Z}{\ensuremath{\mathbf Z}}
\newcommand{\Q}{\ensuremath{\mathbf Q}}
\newcommand{\R}{\ensuremath{\mathbf R}}
\newcommand{\C}{\ensuremath{\mathbf C}}

\newcommand{\bsmallmat}{\left(\begin{smallmatrix}}
\newcommand{\esmallmat}{\end{smallmatrix}\right)}

\DeclareMathOperator{\id}{id}

\newcommand{\bmat}{\left(\begin{matrix}}
\newcommand{\emat}{\end{matrix}\right)}

\DeclareMathOperator{\qqand}{\qquad\text{and}\qquad}

\DeclareMathOperator{\image}{\text{image}}

\DeclareMathOperator{\Log}{Log}

\title{Exponential automorphisms and a problem of Mycielski}
\author{Melvyn B. Nathanson}
\address{Department of Mathematics\\Lehman College (CUNY)\\
Bronx, NY 10468}
\email{melvyn.nathanson@lehman.cuny.edu}
\date{\today}



\thanks{Supported in part by PSC-CUNY Research Award Program grant 63117-00 51.}

\begin{document}

\begin{abstract} An exponential automorphism of $\mathbf{C}$ is a function 
$\alpha: \mathbf{C} \rightarrow \mathbf{C}$  such that  
$\alpha(z_1 + z_2) = \alpha(z_1) + \alpha(z_2)$ and 
$\alpha\left( e^z \right) = e^{\alpha(z)}$ for all $z, z_1, z_2 \in \C$.  
Jan Mycielski  asked if 
$\alpha(\ln 2) = \ln 2$ and if $\alpha(2^{1/k}) = 2^{1/k}$ for $k = 2, 3, 4$
and for all exponential automorphisms $\alpha$.  
These questions are answered modulo a multiple of $2\pi i$ and  a root of unity.    
\end{abstract}

\maketitle

\section{Exponential automorphisms} 

An \emph{automorphism} of the field $K$  
is a function $\sigma: K \rightarrow K$ 
 such that 
\[
\sigma(z_1 + z_2) = \sigma(z_1) + \sigma(z_2)
\]
and
\[
\sigma(z_1 z_2) = \sigma(z_1)+ \sigma(z_2)
\]
for all $z_1, z_2 \in K$, and also 
  \[
 \sigma(1) = 1. 
 \] 

An \emph{exponential automorphism} of the field \C\ 
of complex numbers is a function $\alpha: \mathbf{C} \rightarrow \mathbf{C}$ 
 such that 
\beq                          \label{condition1}
\alpha(z_1 + z_2) = \alpha(z_1) + \alpha(z_2)
\eeq
and
\beq                          \label{condition2}
\alpha\left( e^z \right) = e^{\alpha(z)}
\eeq
for all $z,z_1, z_2 \in \C$. 
The identity function $z \mapsto z$ and complex conjugation $z \mapsto \overline{z}$ 
are examples of automorphisms and exponential automorphisms.

Let $\ln r$ denote the usual real-valued logarithm of the positive real number $r$.
For every positive integer $k$, let $2^{1/k}$ be the positive $k$th root of 2.  
Jan Mycielski~\cite{myci22}
 asked if, for every exponential automorphism $\alpha$, 
 we have $\alpha\left(2^{1/k} \right) = 2^{1/k}$ for $k = 2, 3, 4$ 
 and $\alpha(\ln 2) = \ln 2$. 
He observed that the answers to these questions are ``yes''   
if $\alpha$ is the identity function or if $\alpha$ is complex conjugation.

\bt                          \label{Myc:theorem:Q}    
If  $\alpha$ is an exponential automorphism, then        
\[
 \alpha \left( t \right) = t
\]
for all $t \in \Q$.
\et

\begin{proof}  
We have 
\[
\alpha(1) = \alpha(1+0) = \alpha(1)  + \alpha(0) 
\]
and so
\[
\alpha(0) = 0
\]
and 
\[
\alpha(1) = \alpha\left( e^0 \right) = e^{\alpha(0)} = e^0 = 1.
\]
For all positive integers $m$ and $n$ we have 
\[
m = m\alpha(1) =  \alpha(\underbrace{1 + \cdots + 1 }_{m \text{ summands}}) = 
\alpha(m) = \alpha \left( \underbrace{\frac{m}{n} + \cdots +  \frac{m}{n}  }_{n \text{ summands}} \right) 
= n \alpha \left(  \frac{m}{n}  \right)
\]
and so 
\[
 \alpha(m) = m \qqand \alpha \left(  \frac{m}{n}  \right) =   \frac{m}{n}.
\]
For all $z \in \mathbf{C}$ we have 
\[
0 = \alpha(0) = \alpha(z + (-z)) = \alpha(z) + \alpha(-z)
\]
and so $\alpha(-z) = -\alpha(z)$.  Thus, $\alpha(t) = t$ for all $t \in \mathbf{Q}$. 
This completes the proof. 
\end{proof}

\bt                 \label{Myc:theorem:ln-2}
For every exponential automorphism $\alpha$  
and for every positive rational number $t$ there is an integer 
$n = n(\alpha,t)$ such that 
\[
 \alpha(\ln t) = \ln t + 2\pi n i.
 \]
\et

\begin{proof}
We have 
\[
t = \alpha(t) = \alpha\left( e^{\ln t}\right) = e^{\alpha(\ln t)} 
\]
and so 
\[
 \alpha(\ln t)  = \ln t + 2\pi n i 
\] 
for some integer $n$. 
This completes the proof.
\end{proof} 

Thus, $\alpha(\ln 2) = \ln 2 + 2\pi ni$ for some integer $n$.

\section{Review of logarithms}
The complex logarithm $\log z$ is a  function from  the set $\C^{\times}$ 
of nonzero complex numbers
to the set of sets of complex numbers:   
If $z = re^{i\theta}  \in\C^{\times}$ 
with $r, \theta \in \R$ and $r > 0$, then  
\[
\log z = \left\{\ln r + i\theta + 2\pi n i : n \in \Z  \right\}.  
\]
For all $n \in \Z$ we have 
\[
e^{\ln r + i\theta + 2\pi n i} = e^{\ln r} e^{i\theta} e^{2\pi n i} = r e^{i\theta} = z.
\]
Thus, we can define $e^{\log z} = z$ for all $z \in \C^{\times}$.

For all $\theta \in \R$ there is a unique integer $n$ such that 
$-\pi < \theta + 2\pi n \leq \pi$.
The \emph{principal branch} of the complex logarithm  is the function 
$\Log  : \C^{\times} \rightarrow \C$ defined by 
\[
\Log   z = \ln   r + i\theta \in \log z \qquad \text{with $-\pi < \theta \leq \pi$.}
\]
We have  
\beq                       \label{Myc:logz}   
\log z = \left\{ \Log z  + 2\pi n i : n \in \Z  \right\}  
\eeq
and $\Log r = \ln r$ for all $r > 0$.

Let $z_1 = r_1e^{\theta_1 i}$ and $z_2 = r_2e^{\theta_2 i}$, with
$-\pi < \theta_1, \theta_2 \leq \pi$.  We have 
\[
\Log z_1 = \ln r_1 + i\theta_1 \qqand \Log z_2 = \ln r_2 + i\theta_2  
\]
and 
\[
-2\pi < \theta_1+ \theta_2 \leq 2\pi.
\]
There is a unique integer $\delta \in \{0,1,-1\}$ such that 
\[
-\pi < \theta_1+ \theta_2 +  2\pi \delta   \leq \pi.
\]
Let $ \theta_3 = \theta_1+ \theta_2 +  2\pi \delta$.  
It follows that 
\[
z_1z_2 = r_1r_2e^{(\theta_1+\theta_2)i} = r_1r_2e^{(\theta_3 - 2\pi \delta )i } 
= r_1r_2e^{\theta_3i }
\]
 and so 
\begin{align}                \label{Myc:logz1z2}
\Log z_1 z_2 & = \ln r_1 r_2 +  \theta_3 i    \\
& = \ln r_1 + \ln  r_2 +  (\theta_1 + \theta_2 +  2\pi \delta)i     \nonumber \\
& = \ln r_1 + \theta_1 i +  \ln  r_2 + \theta_2 i + 2\pi\delta i     \nonumber   \\
& = \Log z_1 + \Log z_2 + 2\pi \delta i.        \nonumber 
\end{align}

\section{From exponential automorphisms to automorphisms} 

\bt                                                                  \label{Myc:theorem:automorphism} 
Every exponential automorphism $\alpha$ is an automorphism of \C.
\et

\begin{proof}
By Theorem~\ref{Myc:theorem:Q}, the restriction of  the function 
$\alpha$ to \Q\ is the identity.  
Because $\alpha$ satisfies the additivity relation~\eqref{condition1}, 
it suffices to prove that 
\[
 \alpha(z_1z_2) =  \alpha(z_1) \alpha(z_2)
\]
for all $z_1, z_2 \in \C^{\times}$.  
By
~\eqref{Myc:logz1z2}, there exists $\delta = \delta(z_1, z_2) \in \{0,1,-1\}$ 
such that 
\begin{align*}
\alpha(z_1z_2)  & = \alpha \left( e^{\Log z_1z_2} \right) \\
& =  e^{\alpha \left( \Log z_1z_2 \right) } \\
& = e^{\alpha \left( \Log  z_1 + \Log  z_2 + 2\pi  \delta   i \right)} \\ 
& = e^{\alpha  ( \Log  z_1) + \alpha(\Log  z_2 ) + \alpha(2\pi  \delta   i) } \\
& = e^{\alpha  ( \Log  z_1)} e^{ \alpha(\Log  z_2 ) }   e^{ \alpha(2\pi  \delta   i )} \\ 
& =  \alpha \left( e^{\Log  z_1} \right)  \alpha \left( e^{\Log  z_2} \right)   \alpha \left( e^{2\pi  \delta   i} \right) \\
& = \alpha(z_1) \alpha(z_2) \alpha(1) \\
& = \alpha(z_1) \alpha(z_2).
\end{align*} 
This completes the proof.
\end{proof}

\bt
Let $r$ be a positive rational number and  let $k$ be a positive integer. 
If $r^{1/k}$ is the positive $k$th root of $r$, then 
\[
\alpha\left(r^{1/k} \right) = \zeta r^{1/k}
\]
where $\zeta$ is a $k$th root of unity. 
\et

\begin{proof}
Let  $t = r^{1/k}$.  We have 
\[
r = \alpha(r) = \alpha\left( t^k\right) = \alpha(t)^k 
\]
and so 
\[
\alpha\left(r^{1/k} \right) = \alpha\left( t \right) =\zeta r^{1/k}
\]
for some $k$th root of unity $\zeta$ . 
\end{proof}

\section{An open problem}

\bt
Let $\alpha$ be an exponential automorphism.  If $\alpha(x) \in \R$ for all $x \in \R$  
or if  $\alpha$ is continuous, then the answer to Mycielski's question is ``yes.'' 
\et

\begin{proof}
The following argument uses classical results in algebra (cf. Yale~\cite{yale66}),
whose proofs are in the Appendix.

By Theorem~\ref{Myc:theorem:automorphism}, 
every exponential automorphism $\alpha$ of \C\ is an automorphism of \C. 
By Theorem~\ref{Myc:theorem:R-auto}, 
the only automorphism of \R\ is the identity.  
Thus, if $\alpha$ is an exponential automorphism 
such that $\alpha(x) \in \R$ for all $x \in \R$, 
then $\alpha$ is the identity on \R. 
By Theorem~\ref{Mycielski:theorem:C-auto}, 
the only automorphisms of \C\ whose restriction to  \R\ is the identity 
are the identity and complex conjugation. 
By Theorem~\ref{Mycielski:theorem:C-cont},  
the only continuous automorphisms of \C\ are 
the identity and complex conjugation.  The answer to Mycielski's 
question is ``yes'' for both these automorphisms. 
\end{proof}

It is natural to ask if the only exponential automorphisms  are the identity 
function and complex conjugation.  
Equivalently,  is every exponential automorphism real-valued on \R?  
Is every exponential automorphism continuous?

\appendix
\section{Automorphisms of \R\ and \C}

\bt                   \label{Myc:theorem:R-auto}
The only automorphism of the field \R\ is the identity.
\et

\begin{proof}
Let $\sigma$ be an automorphism of \R. 
We have $\sigma(1) = 1$.  It follows from 
the proof of Theorem~\ref{Myc:theorem:Q} that 
$\sigma(t) = t$ for all $t \in \Q$.

If $x \in \R$ and  $x > 0$, then $x = y^2$ for some nonzero $y \in \R$ 
and so $\sigma(x) = \sigma\left(y^2\right) = \sigma(y)^2 > 0$.
If $x,y \in \R$ and $x< y$, then $y-x>0$ and 
\[
\sigma(y) - \sigma(x) = \sigma(y-x) > 0.
\]
Thus, $x < y$ implies $\sigma(x) < \sigma(y)$ and the function $\sigma$ is order preserving.

If $x \in \R$  and $\sigma(x) \neq x$, then either $x < \sigma(x)$ or $\sigma(x) < x$.
If $x < \sigma(x)$, then there exists $t \in \Q$ such that 
\[
x < t < \sigma(x)
\]
and so 
\[
\sigma(x) < \sigma(t) = t < \sigma(x)
\]
which is absurd. 
If $\sigma(x) < x$, then there exists $t \in \Q$ such that 
\[
\sigma(x) < t < x
\]
and so 
\[
\sigma(x) < t  = \sigma(t) < \sigma(x)
\]
which is also absurd. 
Therefore, $\sigma(x) = x$ for all $x \in \R$.
This completes the proof. 
\end{proof}

\bt        \label{Mycielski:theorem:C-auto}
The only automorphisms $\sigma$ of \C\ such that $\sigma(x) \in \R$ 
for all $x \in \R$ are the identity and complex conjugation. 
\et

\begin{proof}
Let $\sigma$ be an automorphism of \C. 
If $\sigma(x) \in \R$ for all $x \in \R$, then the restriction of $\sigma$ to \R\ 
is an automorphism of \R\ and so $\sigma$ is the identity on \R.  
We have 
 $\sigma(i)^2 = \sigma(i^2) = \sigma(-1)  = -1$ and so $\sigma(i) = i \quad\text{or} \quad \sigma(i) = -i$.
Let  $z = x+yi \in \C$ with $x,y \in \R$.  If $\sigma(i) = i$, then 
\[
\sigma(x+yi) = \sigma(x)+\sigma(y) \sigma(i) = x+yi = z 
\]
and $\sigma$ is the identity on \C.  
If  $\sigma(i) = -i$, then 
\[
\sigma(x+yi) = \sigma(x)+  \sigma(y) \sigma(i)= x+y(-i) = x-yi = \overline{z}
\]
and $\sigma$ is complex conjugation on \C. 
This completes the proof. 
\end{proof}

\bt                  \label{Mycielski:theorem:C-cont}
The only continuous automorphisms of \C\ are the identity  
and complex conjugation. 
\et

\begin{proof}
Let $\sigma$ be a continuous automorphism of \C.  
We have  $\sigma(t) = t$ for all $t \in \Q$ and $\sigma(i) =  \pm i$.  
If $\sigma(i) = i$, then $\sigma(t_1+it_2) = t_1+it_2$ for all Gaussian numbers $t_1+it_2 \in \Q(i)$. 
If $\sigma(i) = -i$, then $\sigma(t_1+it_2) = t_1 - it_2$ for all Gaussian numbers $t_1+it_2 \in \Q(i)$. 

The Gaussian numbers are dense in \C.  For every complex nunber $z$ 
there is a sequence of Gaussian numbers $(w_n)_{n=1}^{\infty}$ such that 
\[
z = \lim_{n\rightarrow \infty} w_n.
\]
If $\sigma$ is a continuous automorphism of \C, then 
\begin{align*}
\sigma(z) & = \sigma\left( \lim_{n\rightarrow \infty} w_n \right) = \lim_{n\rightarrow \infty}  \sigma\left( w_n \right) \\
& = \begin{cases}
\lim_{n\rightarrow \infty} w_n  & \text{if $\sigma(i) = i$}\\ 
\lim_{n\rightarrow \infty} \overline{w_n}  & \text{if $\sigma(i) = -i$}
\end{cases} \\
& = \begin{cases}
 z & \text{if $\sigma(i) = i$}\\ 
 \overline{z} & \text{if $\sigma(i) = -i$}
\end{cases} 
\end{align*}
and so $\sigma$ is either the identity function or complex conjugation. 
This completes the proof. 
\end{proof}


\begin{thebibliography}{99}
\bibitem{myci22}
J. Mycielski, Problem 12301,  Amer. Math. Monthly 129 (2022), 186.

\bibitem{yale66}
P. B. Yale, Automorphisms of the complex numbers, 
Math. Mag. 39 (1966), 135--141.

\end{thebibliography}
\end{document}